\documentclass[onecolumn]{IEEEtran}
\usepackage[colorlinks=true, pdfstartview=FitV, linkcolor=black, citecolor= black, urlcolor= black]{hyperref}
\usepackage{graphicx}    
\usepackage{amssymb}
\usepackage{subfig}
\usepackage{natbib}
\usepackage{multirow}
\usepackage{bm}
\usepackage{pmat}
\usepackage{amsmath,amsthm,amsfonts,amscd} 
\usepackage{fouridx}
\usepackage{eucal} 	 	
\usepackage{verbatim}

\newtheorem{thm}{Theorem}[section]
\newtheorem{lem}{Lemma}[section]
\newtheorem{cor}{Corollary}[section]

\theoremstyle{remark}

\theoremstyle{definition}
\newtheorem{not1}{Note}[section]

\allowdisplaybreaks[3]
\newcommand{\ud}{\mathrm{d}}

\newcommand{\xv}{\bm{x}}

\newcommand{\zv}{\bm{z}}

\begin{document}

%

\title{Stabilizing Controllers for Multi-Input, Singular Control Gain Systems$^{\dagger} $ \thanks{${\dagger} $ A version of this paper was presented at the IEEE Conference in Decision and Control, 2012 in Maui, Hawaii}}

\author{Srikant Sukumar$^{a}$ $^{\ddagger} $ \thanks{$\ddagger $Corresponding Author,  Email: srikant@sc.iitb.ac.in} and Maruthi R. Akella$^{b}$ $^{\S}$ \thanks{$\S$ Email: makella@mail.utexas.edu}\\
\vspace{6pt}
$^{a}$ {\em Systems and Control Engineering, Indian Institute of Technology, Bombay,Mumbai, India}\\
$^{b}${\em Aerospace Engineering and Engineering Mechanics, The University of Texas at Austin, Austin, Texas, USA}\\
\vspace{6pt} 
}


%


\maketitle

\begin{abstract}
This paper proposes a new methodology for design of a stabilizing control law for multi-input linear systems with time-varying, singular gains on the control. The results presented here assume the control gain to satisfy persistence of excitation which is a necessary condition for existence of stabilizing controllers in the presence of unstable drift. This work involves a novel persistence filter construction and provides a significant extension to the authors' previous result on stabilization of single-input linear systems with time-varying singular gains. An application to underactuated spacecraft stabilization is shown which illustrates the interesting features of the time-varying control design in stabilization of nonlinear dynamical systems. Finally, the development of an observer counterpart of these results is presented in the presence of multiple-outputs subject to singular measurement gains.
\end{abstract}


%

\section{Introduction} \label{intro}

The interest in designing feedback controllers for systems with singular, state or time-varying control gain stems from several representative real world applications. For any given $n\times m$ control gain matrix $B(t)$ with $m \le n$, ``singularity'' in the context of this paper refers to the outer-product $B(t)B^T(t)$ having less than minimal rank (i.e., $< m$) at various time instants or even possibly over a finite time-period. An application pertaining to aerospace engineering is the spacecraft attitude stabilization problem using only magnetic torquers as actuators for Low Earth Orbit (LEO) satellites. Magnetic actuation is also of interest in space based interferometry missions such as Terrestrial Planet Finder (TPF)~\footnote{\url{http://planetquest.jpl.nasa.gov/TPF-I/tpf-I_index.cfm}} and the MicroArcsecond X-ray Imaging Mission (MAXIM)~\footnote{\url{http://maxim.gsfc.nasa.gov/mission/mission.html}}. Available solutions to the attitude control problem with magnetic actuation involve linearization~\citep{alfriend76}, periodicity assumption on magnetic field vector~\citep{psiaki01} or control design based on the time-averaged dynamics~\citep{astolfi04}.

Another class of differential equations possessing singular, state or time dependent control gains are nonholonomic chain integrators. 
\citet{kolmanovsky02} and \citet{murray94} provide an extensive survey of the work done in chain form systems. 
Further, stabilizing control design for these classes of systems using notions of persistence of excitation of the control gain have been illustrated by \citet{loria02}. However, no stabilizing control designs are available for more general classes of chain form systems with unstable drift to best of the authors' knowledge.

Examples of state and/or time-varying gains also arise in the area of biomechanical engineering and pertains to Functional Electrical Stimulation (FES). Functional Electrical Stimulation can be concisely defined as the process of applying electrical pulses to nerve fibers resulting in muscle activation. 
Dynamical models of the knee-shank system with time-dependent control gains have been proposed in FES literature for example by~\citet{durfee93}.

In addition to dynamics where state and time-varying control or measurement gains arise naturally, it is also possible to conceive applications wherein gains are artificially introduced to suit specific design requirements. Time-varying gains on the actuator can be utilized to schedule actuator operation, e.g. to allow for intermittent actuator operation due to mission or hardware constraints. Specific examples include: (i) flow control of scramjet engines,
(ii) control of underactuated systems using actuator re-orientation, (iii) simultaneous sensing and actuation with reversible transducers~\citep{srikant12j}.

In view of the above motivating examples and the fact that there remain several unanswered theoretical questions~\citep{loria05}, it is of interest to look at design of controllers for systems with singular, state or time-varying scaling on the control.  The problem of stabilizing single-input dynamics,
\begin{equation} \label{sys0}
    \dot{\xv} = A\xv + g(t)Bu, \quad u(t) \in \mathbb{R}, \, \xv(t_0) = \xv_0 \in \mathbb{R}^n, t_0 \in \mathbb{R}
\end{equation}
with arbitrary unstable drift matrix $A$  and singular gain $g(t)$ was successfully resolved in \citet{srikant09}. In this article we seek to extend the single-input results to multiple input linear dynamics with singular control gains. Specifically, we seek to stabilize dynamics in the form,
\begin{equation} \label{midyn}
	\dot{\xv} = A\xv + BG(t)\bm{u}
\end{equation}
with $\xv \, \in \mathbb{R}^n$, $\bm{u} \in \mathbb{R}^m$ with $m \leq n$, $A \in \mathbb{R}^{n\times n}$, $B \in \mathbb{R}^{n \times m}$ and $G(\cdot):\mathbb{R}^{\geq 0} \rightarrow \mathbb{R}^{m \times m}$ of the form,
\begin{equation} \label{G}
G(t) = \textrm{diag}([g_1(t),\ldots,g_m(t)]^T)
\end{equation}	

Previous attempts at this problem can be attributed to classical works by \citet{morgan77,morgan77b}, \citet{sondhi76} and \citet{krei77}. More recently the stabilization of singular-gain control systems has been studied in-depth and solutions suggested by \citet{chaillet07}. However these have been restricted to the no-drift or neutrally stable drift cases only. 

A notable recent contribution is the work by \citet{weiss12} where the stabilization problem of Linear Time Varying (LTV) systems with singular control gain is studied using a forward Riccatti based control law. The authors attack the basic problem in classical control design for such systems that require solving a backward in time Riccatti equation therefore assuming knowledge of the singular control gain for all time. This is typically not the case and therefore implementation of the classical results is not feasible. The authors demonstrate that for LTV systems with a special structure, (such as the closed-loop dynamics being symmetric or commuting with its integral) a forward Riccatti equation based control law can result in stabilization. A scalar version provided by the authors is similar to \citep{srikant09}. However for general LTV multiple input dynamics, solutions to the singular control gain problem do not exist. In this work we look at feedback design for stabilizing multi-input dynamics with singular control gains for the special case where the control gain matrix is diagonal. The overall closed-loop system however does not satisfy any symmetry or commutativity condition. We also exhibit possible extension of the results presented to the nonlinear spacecraft attitude stabilization problem.

Section~\ref{PEstab} constitutes the main results of this paper and addresses the problem of stabilizing multi-input dynamics with diagonal, matrix time-dependent gains on the control. Section~\ref{mobs} outlines the multi-output observer design counterpart for systems with singular measurement gains. Section~\ref{sim} illustrates application of the linear systems result to a practical nonlinear dynamical system. Attitude and angular velocity stabilization of an axi-symmetric spacecraft with only two independent actuators is considered in this section. The conclusions of this work are summarized in Section~\ref{conc}. 

Classical definitions of Persistence of Excitation (PE) and Exponential Stability (E.S.) as in \citet[p. 24-25, 72]{sastry} are referred throughout this work.

\section{Persistent Filters for Stabilization of Multi-input Linear Systems} \label{PEstab}

\subsection{Canonical Transformations}

%

A specific canonical form for multi-input, multi-output dynamics is employed and  referred to the work by \citet{anderson67}. 
They provide construction of a similarity transformation $T$ using elements of the controllability matrix which yields the canonical system $\hat{A} = TAT^{-1}$ and $\hat{B} = TB$ defined as follows.
\begin{equation} \label{Ahat}
	\hat{A} = \begin{bmatrix}
							\hat{A}_{p,p} & \bm{0} & \bm{0} & \ldots & \bm{0}  \\
							\hat{A}_{p-1,p} & \hat{A}_{p-1,p-1} & \bm{0} & \ldots & \bm{0} \\
							\vdots & \vdots & \ddots & \vdots & \vdots \\
							\hat{A}_{2,p} & \hat{A}_{2,p-1} & \ldots & \hat{A}_{2,2} & \bm{0} \\
							\hat{A}_{1,p} & \hat{A}_{1,p-1} & \ldots & \hat{A}_{1,2} & \hat{A}_{1,1}
						\end{bmatrix} 
\end{equation}						
$\hat{A}$ clearly has a lower triangular block structure with diagonal blocks defined as,
\begin{equation} \label{Ajj}
	\hat{A}_{j,j} = \begin{bmatrix}
										0 & 1 & 0 & \ldots & 0 \\
										0 & 0 & 1 & \ldots & 0 \\
										\vdots & \vdots & \vdots & \ddots & \vdots \\
										0 & \ldots & 0 & 0 & 1 \\
										-\alpha_{j,r_j} & -\alpha_{j,r_j-1}	& \ldots & -\alpha_{j,2} & -\alpha_{j,1}
									\end{bmatrix}
\end{equation}									
which is reminiscent of the controller canonical form for single-input systems utilized for control design in \citep{srikant09}. The lower triangular block matrices in $\hat{A}$ turn out to have the following structure,
\begin{equation} \label{Akj}
	\hat{A}_{k,j} = \begin{bmatrix}
										\beta_{k,j,1} & 0 & \ldots & 0 \\	
										\vdots & \vdots & \vdots & \vdots\\				
										\beta_{k,j,r_k} & 0 & \ldots & 0 
									\end{bmatrix}, \quad \forall k < j
\end{equation}
for some constants $\beta_{k,j,l} \, \text{for} \, l = 1,2 \ldots, r_k$ which signify the nature of the coupling between the states in each block. The lower-triangular structure implies a unidirectional coupling and the first block is completely decoupled from the other blocks. Further, it is interesting to note that the evolution of the states in each block depends only on the first state in each of the previous blocks due to the structure of $\hat{A}_{k,j}$ in~(\ref{Akj}).

The structure of the control scaling matrix $\hat{B}$ for the transformed system comes out to be,
\begin{equation} \label{Bhat}
\begin{tabular}{lc|c}
 & $p$ & $m-p$ \\
$\hat{B} \; = $ & $\left[ \bm{e}_n \; \bm{e}_{n-r_1} \; \bm{e}_{n-r_1-r_2} \; \ldots \; \bm{e}_{r_p}\right.$ & $\left. \times \;\; \times\;\; \right]$
\end{tabular}
\end{equation}	
where $\bm{e}_k$ is a unit vector with zeros for all elements except the $k^{th}$ component which is unity and each of the $\times $'s represent possible non-zero vectors. It is evident from the partitioned matrix $\hat{B}$ that the controls corresponding to the first $p$ (also the number of blocks in $\hat{A}$) columns are sufficient to control the system while the rest $(m-p)$ controls are redundant and will henceforth be set to zero. Finally, the transformed dynamics is given by,

\begin{equation} \label{mizdyn}
	\dot{\bm{z}} = \hat{A}\bm{z} + \hat{B}G(t)\bm{u}
\end{equation}

\subsection{Stabilization of Multi-Input Dynamics}

The results in the preceding section allow control design to be based upon the transformed dynamics~(\ref{mizdyn}). The persistence filter for the multi-input dynamics is a set of linear dynamics identical to the single-input case and corresponding to number of diagonal blocks in $\hat{A}$, i.e.,
\begin{equation} \label{miR}
	\dot{R}_j = -\lambda_j R_j + g^k_j(t), \quad j=1,2,\ldots,p
\end{equation}
with $R_j(0) > 0$, $k = \textrm{max}\{2,2^{\lceil \log_2 r \rceil}\}$, $r = \text{max} \{r_1,r_2,\ldots,r_p\}$ and $\lambda_j > 0$. Since each individual persistence filter is identical to the single-input case, it can be inferred as before from Lemma~4 in \citep{srikant09} and PE of $g_j(t)$ the existence of positive lower and upper bounds on each $R_j(t)$ for all $t\geq 0$. These will henceforth be denoted as, $R_{j,min}$ and $R_{j,max}$ respectively. Persistence filters are required corresponding to the gains on only the last $p$ controls because the rest of the controls are redundant as evidenced by the structure of $\hat{B}$ in Eq.~(\ref{Bhat}). These redundant controls are set to zero and do not play any further role in the system stabilization. In a physical sense, these controls can be eliminated in the design phase itself to avoid over-actuation of the system.

The column vector $\bm{z}$ is partitioned based upon the number of block-diagonal matrices $\hat{A}_{j,j}$ (with $j$ being the assigned block number) as, 
\begin{equation} \label{zpart}
	\bm{z} = [ \fourIdx{p}{}{T}{}{\bm{z}},\;\fourIdx{p-1}{}{T}{}{\bm{z}},\;\ldots,\;\fourIdx{1}{}{T}{}{\bm{z}}]^T 
\end{equation}
where $\fourIdx{j}{}{}{}{\bm{z}} = [\fourIdx{j}{}{}{1}{z},\,\fourIdx{j}{}{}{2}{z},\, \ldots,\,\fourIdx{j}{}{}{r_j}{z}]^T \in \mathbb{R}^{r_j}$ has additional indexing $j$ on components corresponding to the block they belong to. Similar notation will be followed in defining the augmented states $\bm{\Omega}$ which are partitioned as,
\begin{equation} \label{ompart}
	\bm{\Omega} = [\fourIdx{p}{}{T}{}{\bm{\Omega}},\;\fourIdx{p-1}{}{T}{}{\bm{\Omega}},\;\ldots,\;\fourIdx{1}{}{T}{}{\bm{\Omega}}]^T 
\end{equation}

\subsection{Construction of Augmented States}
The construction of augmented states from $\bm{z}$ is an important intermediate step in the control design process and is critical for the Lyapunov analysis of the closed-loop dynamics. The augmented states are defined with respect to each block in matrix $\hat{A}$. Starting with the dynamics of the $p$-block in Eq.~(\ref{mizdyn}),
\begin{equation} \label{mizdyn1}
	\fourIdx{p}{}{}{}{\dot{\bm{z}}} = \hat{A}_{p,p} \fourIdx{p}{}{}{}{\bm{z}} + [0,\;0,\;\ldots,\;g_{p}(t)u_{p}]^T
\end{equation}																																																								where, $u_{p+1}=u_{p+2}=\cdots=u_m=0$ and $\hat{A}_{p,p}$ is defined	in~(\ref{Ajj}). The dynamics are in the standard controller canonical form for single-input systems. The augmented state definition is therefore given by~\citep{srikant09},
\begin{eqnarray} \label{Om1}
	  \fourIdx{p}{}{}{1}{\Omega} &\triangleq& \fourIdx{p}{}{}{1}{z} \nonumber \\
    \fourIdx{p}{}{}{i+1}{\Omega} &\triangleq& \frac{\ud}{\ud t} \fourIdx{p}{}{}{i}{\Omega} + \frac{g^k}{2R} \fourIdx{p}{}{}{i}{\Omega}, \quad i = 1,2,\dots,r_p-1
\end{eqnarray}

Moving forward, the $(p-1)$-block has dynamics given by,
\begin{eqnarray} \label{mizdyn2}
  \fourIdx{p-1}{}{}{}{\dot{\bm{z}}} &=& \hat{A}_{p-1,p-1} \fourIdx{p-1}{}{}{}{\bm{z}} + [\beta_{p-1,p,1},\;\;\ldots,\;\beta_{p-1,p,r_{p-1}} ]^T\fourIdx{p}{}{}{1}{\Omega} \nonumber \\
	& & + [0,\;\ldots,\;g_{p-1}(t)u_{p-1}]^T
\end{eqnarray}		
which as mentioned earlier indicates coupling with the previous block. This coupling motivates a slightly different choice of augmented states for this block as compared to Eq.~(\ref{Om1}),
\begin{eqnarray} \label{Om2}
		\fourIdx{p-1}{}{}{1}{\Omega} &\triangleq& \fourIdx{p-1}{}{}{1}{z} \nonumber \\
    \fourIdx{p-1}{}{}{i+1}{\Omega} &\triangleq& \frac{\ud}{\ud t} \fourIdx{p-1}{}{}{i}{\Omega} + \frac{g^k}{2R} \fourIdx{p-1}{}{}{i}{\Omega} - \beta_{p-1,p,i}\fourIdx{p}{}{}{1}{\Omega}
\end{eqnarray}
for $i = [1,2,\dots,r_{p-1}-1$.

Following the same pattern as above, the augmented states for the $(p-j)$-block ($j < p$) are defined as,
\begin{eqnarray} \label{Omj}
		\fourIdx{p-j}{}{}{1}{\Omega} &\triangleq & \fourIdx{p-j}{}{}{1}{z} \nonumber \\
    \fourIdx{p-j}{}{}{i+1}{\Omega} &\triangleq & \frac{\ud}{\ud t} \fourIdx{p-j}{}{}{i}{\Omega} + \frac{g^k}{2R} \fourIdx{p-j}{}{}{i}{\Omega} - \sum_{s=p-j+1}^{p}\beta_{p-j,s,i}\fourIdx{s}{}{}{1}{\Omega} \nonumber \\
    & &
\end{eqnarray}
for $i = 1,2,\dots,r_{p-j}-1$.

\subsection{Feedback Law Design}
The stabilization result involving the persistence filter and augmented state definitions is stated as a theorem below followed by a Lyapunov analysis based proof of the same.

\begin{thm} \label{miPEthm1}

Consider the linear, multi-input dynamics in Eq.~(\ref{midyn}) under the assumptions that $(A,B)$ is a controllable pair and the component functions of $G(t)$ are $\mathcal{C}^{n-1}$, bounded with bounded derivatives up to order $(n-1)$ and satisfy the PE condition. Further, let $\bm{z} = T\xv$ where $T$ is the transformation to the canonical form defined in Eqs.~(\ref{Ahat})-(\ref{Bhat}) and the augmented states defined in Eqs.~(\ref{Om1}) and (\ref{Omj}). Then the following control law,
\begin{equation} \label{mipcont}
        u_j = -\frac{1}{g_j}\sum_{i=1}^{r_j}\alpha_{j,r_j-i+1}(\fourIdx{j}{}{}{i}{\Omega}-\fourIdx{j}{}{}{i}{z})- \frac{1}{g_j}\frac{\ud}{\ud t} (\fourIdx{j}{}{}{r_j}{\Omega}-\fourIdx{j}{}{}{r_j}{z}) - \frac{g^{k-1}_j}{2R_j}\fourIdx{j}{}{}{r_j}{\Omega}	
\end{equation}
for $j=1,2,\ldots,p$ and $u_j = 0$ for $j>p$ with the persistence filters defined in Eq.~(\ref{miR}) with $R_j(0) > 0$, $k = \textrm{max}\{2,2^{\lceil \log_2 r \rceil}\}$, $r = \text{max} \{r_1,r_2,\ldots,r_p\}$ and $\lambda_j > 0$ guarantees exponential convergence of $\xv$ to the origin subject to following inequalities on $\lambda_j,\,j=1,2,\ldots,p$,
\begin{eqnarray} \label{gamjineq}
	\gamma_1 &>& \sum_{l=1}^{p-1} \|\hat{A}_{1,1+l}\| \nonumber \\
	\gamma_j &>& \sum_{l=1}^{p-j} \|\hat{A}_{j,j+l}\|+ 2\sum_{l=1}^{j-1} \|\hat{A}_{j-l,j}\| , \: j = 2,\ldots p-1 \nonumber \\
	\gamma_p &>& \sum_{l=1}^{p-1} \|\hat{A}_{p-l,p}\| 
\end{eqnarray}
where, $\gamma_j$ for $j=1,2,\ldots,p$ is defined to be,
\begin{eqnarray} \label{gamj}
	\gamma_j &\triangleq & \lambda_j - \max_{l\in[2,r_j-1]}\bigg\{ (1+|\alpha_{j,r_j}|),(2+|\alpha_{j,r_j-l+1}|), (1-\alpha_{j,1} \nonumber \\
					& & +\sum_{i=1}^{r_j-1}|\alpha_{j,r_j-i+1}|)\bigg\} \triangleq  \lambda_j - \lambda_j^*
\end{eqnarray}
	
The rate of exponential convergence can be made arbitrarily large by appropriate choice of persistence filter gains and the following expression provides an estimate for the convergence rate.
\begin{equation} \label{sig}
	\sigma = \min\left\{ \frac{\gamma_p}{2},\frac{\gamma_{p-j}- \sum_{l=1}^{p-j} \|\hat{A}_{j,j+l}\|}{2},\gamma_1 - \sum_{l=1}^{p-1} \|\hat{A}_{1,1+l}\|\right \}
\end{equation}

\end{thm}    

\begin{not1}
	The structure of the controller is identical to the single-input feedback law. Remark~9 in \cite{srikant09} therefore ensures that the division by singular gain $g_j$ in the control law above is only \emph{symbolic}. The choice of $k$ specified after Eq.~(\ref{miR}) ensures that each term above is scaled by $g_j^\gamma,\, \gamma \geq 1$ guaranteeing cancellation of $g_j$ in the numerator and denominator.
\end{not1}
	
\begin{proof}
Corresponding to each block in $\hat{A}$, energy functionals are defined which are then combined with appropriate scaling to arrive at a candidate Lyapunov function for the entire dynamics~(\ref{mizdyn}). The constituent energy functional for arbitrary block $j$ where $j = 1,2,\ldots,p$ is,
\begin{equation} \label{miVo}
	V_j^o = R_j\sum_{i=1}^{r_j} \fourIdx{j}{}{2}{i}{\Omega} 
\end{equation}
which has the following derivative accounting for the persistence filter dynamics~(\ref{miR}),
\begin{equation} \label{miVod}
	\dot{V}_j^o = -\lambda_j R_j \sum_{i=1}^{r_j}  \fourIdx{j}{}{2}{i}{\Omega} + 2R_j \sum_{i=1}^{r_j} \fourIdx{j}{}{}{i}{\Omega}\left( \frac{\ud}{\ud t} \fourIdx{j}{}{}{i}{\Omega} + \frac{g^k_j}{2R_j} \fourIdx{j}{}{}{i}{\Omega} \right)
\end{equation}

Focusing now on each individual block, the mixed term in Eq.~(\ref{miVod}) for the $p$-block is computed. From the augmented state definitions for this block in Eq.~(\ref{Om1}) it can be shown that,
\begin{align} \label{miVodp1}
	& \sum_{i=1}^{r_p} \fourIdx{p}{}{}{i}{\Omega}\left( \frac{\ud}{\ud t} \fourIdx{p}{}{}{i}{\Omega} + \frac{g^k_p}{2R_p} \fourIdx{p}{}{}{i}{\Omega} \right) = \sum_{i=1}^{r_p-1} \fourIdx{p}{}{}{i}{\Omega} \fourIdx{p}{}{}{i+1}{\Omega} +\nonumber \\
	&	\fourIdx{p}{}{}{r_p}{\Omega} \left( \frac{\ud}{\ud t}  \fourIdx{p}{}{}{r_p}{\Omega} + \frac{g^k_p}{2R_p} \fourIdx{p}{}{}{r_p}{\Omega} \right) & &
\end{align}	
where the last term on the right hand side above can be evaluated from dynamics~(\ref{mizdyn1}) as,
\begin{align} \label{miVodp2}
& \left( \frac{\ud}{\ud t} \fourIdx{p}{}{}{r_p}{\Omega} + \frac{g^k_p}{2R_p} \fourIdx{p}{}{}{r_p}{\Omega} \right)	= \bigg(-\sum_{i=1}^{r_p}\alpha_{p,r_p-i+1}\fourIdx{p}{}{}{i}{z} + g_p u_p + \nonumber \\
& \frac{\ud}{\ud t} (\fourIdx{p}{}{}{r_p}{\Omega}-\fourIdx{p}{}{}{r_p}{z}) +  \frac{g^k_p}{2R_p} \fourIdx{p}{}{}{r_p}{\Omega} \bigg) & & 
\end{align}
Substituting for control $u_p$ in the above expression from Eq.~(\ref{mipcont}) yields,
\begin{equation} \label{miVodp3}
\left( \frac{\ud}{\ud t} \fourIdx{p}{}{}{r_p}{\Omega} + \frac{g^k_p}{2R_p} \fourIdx{p}{}{}{r_p}{\Omega} \right)	= -\sum_{i=1}^{r_p}\alpha_{p,r_p-i+1}\fourIdx{p}{}{}{i}{\Omega} 
\end{equation}
Finally substituting back into Eq.~(\ref{miVodp1}) and then in (\ref{miVod}), the directional derivative $\dot{V}_p^o$ for the $p$-block after applying inequality $2|\omega_1| |\omega_2| \leq \omega_1^2 + \omega_2^2$ can be obtained to be,
\begin{align} \label{miVodp5}
&	\dot{V}_p^o \leq  -\Big[(\lambda_p-1-|\alpha_{p,r_p}|)\, \fourIdx{p}{}{2}{1}{\Omega}\nonumber \\ &+\sum_{i=2}^{r_p-1} (\lambda_p-2-|\alpha_{p,r_p-i+1}|)\, \fourIdx{p}{}{2}{i}{\Omega}  \nonumber \\
&  + (\lambda_p-1+\alpha_{p,1}-\sum_{i=1}^{r_p-1}|\alpha_{p,r_p-i+1}|)\, \fourIdx{p}{}{2}{r_p}{\Omega}\Big]R_p
\end{align}
Now assuming that $\lambda_p$ is chosen large enough to render the bracketed terms positive, according to Eq.~(\ref{gamj}), yields,
\begin{equation} \label{miVodp5.1}
	\dot{V}_p^o \leq -\gamma_p V_p^o = -\gamma_p R_p\sum_{i=1}^{r_p} \fourIdx{p}{}{2}{i}{\Omega} 
\end{equation}

The directional derivative is now computed as shown in Eq.~(\ref{miVod}) for any arbitrary $j < p$. Using the augmented state definitions~(\ref{Omj}), the last term in Eq.~(\ref{miVod}) can be evaluated as,
\begin{align} \label{miVodp6}
	&	\sum_{i=1}^{r_j} \fourIdx{j}{}{}{i}{\Omega}\left( \frac{\ud}{\ud t} \fourIdx{j}{}{}{i}{\Omega} + \frac{g^k_j}{2R_j} \fourIdx{j}{}{}{i}{\Omega} \right) = \sum_{i=1}^{r_j-1} \fourIdx{j}{}{}{i}{\Omega} \fourIdx{j}{}{}{i+1}{\Omega} + \fourIdx{j}{}{}{r_j}{\Omega} \left( \frac{\ud}{\ud t} \fourIdx{j}{}{}{r_j}{\Omega}\right.  \nonumber \\
 & +\left. \frac{g^k_j}{2R_j} \fourIdx{j}{}{}{r_j}{\Omega} \right) + \sum_{i=1}^{r_j-1} \fourIdx{j}{}{}{i}{\Omega} \sum_{s=j+1}^{p}\beta_{j,s,i}\fourIdx{s}{}{}{1}{\Omega}
\end{align}
where the term corresponding to the last augmented state in above equation can be computed as before as,
\begin{align} \label{miVodp7}
& \left( \frac{\ud}{\ud t} \fourIdx{j}{}{}{r_j}{\Omega} + \frac{g^k_j}{2R_j} \fourIdx{j}{}{}{r_j}{\Omega} \right)	= \Big(-\sum_{i=1}^{r_j}\alpha_{j,r_j-i+1}\fourIdx{j}{}{}{i}{z} + g_j u_j  \nonumber \\
&+ \sum_{s=j+1}^{p}\beta_{j,s,r_j}\fourIdx{s}{}{}{1}{\Omega} +\frac{\ud}{\ud t} (\fourIdx{j}{}{}{r_j}{\Omega}-\fourIdx{j}{}{}{r_j}{z}) +  \frac{g^k_j}{2R_j} \fourIdx{j}{}{}{r_j}{\Omega} \Big)
\end{align}
Substituting for control $u_j$ in the above expression from Eq.~(\ref{mipcont}) yields,
\begin{equation} \label{miVodp8}
\left( \frac{\ud}{\ud t} \fourIdx{j}{}{}{r_j}{\Omega} + \frac{g^k_j}{2R_j} \fourIdx{j}{}{}{r_j}{\Omega} \right)	= -\sum_{i=1}^{r_j}\alpha_{j,r_j-i+1}\fourIdx{j}{}{}{i}{\Omega} + \sum_{s=j+1}^{p}\beta_{j,s,r_j}\fourIdx{s}{}{}{1}{\Omega}
\end{equation}
Cross terms begin to appear in Eqs.~(\ref{miVodp6}) and (\ref{miVodp8}) as a consequence of $\hat{A}_{k,j}$ (for $k<j$) matrices being non-zero which indicate coupling with the states in the previous block. However, as stated before it can be verified that the coupling is unidirectional and involves only the first state in each block. Combining Eqs.~(\ref{miVodp8}) and (\ref{miVodp6}) and substituting the result back in~(\ref{miVod}) yields,
\begin{align} \label{miVodp9}
&	\dot{V}_j^o = -\lambda_j R_j \sum_{i=1}^{r_j}  \fourIdx{j}{}{2}{i}{\Omega} + 2R_j \Big[\sum_{i=1}^{r_j-1} \fourIdx{j}{}{}{i}{\Omega} \fourIdx{j}{}{}{i+1}{\Omega} \nonumber \\
& + \sum_{i=1}^{r_j-1} \fourIdx{j}{}{}{i}{\Omega} \sum_{s=j+1}^{p}\beta_{j,s,i}\fourIdx{s}{}{}{1}{\Omega}+ \fourIdx{j}{}{}{r_j}{\Omega}\left(-\sum_{i=1}^{r_j}\alpha_{j,r_j-i+1}\fourIdx{j}{}{}{i}{\Omega} \right. \nonumber \\
&\left. + \sum_{s=j+1}^{p}\beta_{j,s,r_j}\fourIdx{s}{}{}{1}{\Omega}\right) \Big]
\end{align}
The mixed terms corresponding to the $j$-block can be dominated using the negative quadratic terms in the above equation as before and simplified to yield,
\begin{align} \label{miVodp11}
	& \dot{V}_j^o \leq  -R_j\Big[(\lambda_j-1-|\alpha_{j,r_j}|)\, \fourIdx{j}{}{2}{1}{\Omega} + \sum_{i=2}^{r_j-1} (\lambda_j-2 \nonumber \\
	& -|\alpha_{j,r_j-i+1}|)\, \fourIdx{j}{}{2}{i}{\Omega} + (\lambda_j-1+\alpha_{j,1}-\sum_{i=1}^{r_j-1}|\alpha_{j,r_j-i+1}|)\, \fourIdx{j}{}{2}{j}{\Omega}\Big] \nonumber \\
	& + 2R_j\left[ \fourIdx{j}{}{T}{}{\bm{\Omega}} \left(\hat{A}_{j,j+1}\fourIdx{j+1}{}{}{}{\bm{\Omega}}+\ldots + \hat{A}_{j,p}\fourIdx{p}{}{}{}{\bm{\Omega}} \right)  \right]
\end{align}		
Choosing as before $\gamma_j > 0$ as defined in Eq.~(\ref{gamj}) simplifies~(\ref{miVodp11}) to,
\begin{align} \label{miVodp11.1}
&	\dot{V}_j^o \leq  -\gamma_j R_j\sum_{i=1}^{r_j} \fourIdx{j}{}{2}{i}{\Omega} + 2R_j\left[ \fourIdx{j}{}{T}{}{\bm{\Omega}} \left(\hat{A}_{j,j+1}\fourIdx{j+1}{}{}{}{\bm{\Omega}}+\ldots + \hat{A}_{j,p}\fourIdx{p}{}{}{}{\bm{\Omega}} \right)  \right]
\end{align}	

The following energy-like function is defined for combining the $p$ and $p-1$ blocks,
\begin{equation} \label{V12}
	V_{[p,p-1]} \triangleq V_{p-1}^o + 2\frac{R_{p-1,max}}{R_{p,min}} V_p^o
\end{equation}	
The direction derivative of $V_{[p,p-1]}$ can be computed based on Eqs.~(\ref{miVodp5.1}) and (\ref{miVodp11.1}) to be,
\begin{align} \label{V12d}
&	\dot{V}_{[p,p-1]} \leq  -2\frac{R_{p-1,max}}{R_{p,min}}\gamma_p R_p\sum_{i=1}^{r_p} \fourIdx{p}{}{2}{i}{\Omega} -\gamma_{p-1} R_{p-1}\sum_{i=1}^{r_{p-1}} \fourIdx{p-1}{}{2}{i}{\Omega} \nonumber \\
& + 2R_{p-1}\fourIdx{p-1}{}{T}{}{\bm{\Omega}}\hat{A}_{p-1,p}\fourIdx{p}{}{}{}{\bm{\Omega}} \nonumber \\
&\leq -\frac{R_{p-1,max}}{R_{p,min}}\gamma_p R_p\sum_{i=1}^{r_p} \fourIdx{p}{}{2}{i}{\Omega} -\gamma_p R_{p-1}\sum_{i=1}^{r_p} \fourIdx{p}{}{2}{i}{\Omega}\nonumber \\
&  -\gamma_{p-1} R_{p-1}\sum_{i=1}^{r_{p-1}} \fourIdx{p-1}{}{2}{i}{\Omega} +  R_{p-1} \| \hat{A}_{p-1,p} \| \left( \| \fourIdx{p-1}{}{}{}{\bm{\Omega}} \|^2 + \| \fourIdx{p}{}{}{}{\bm{\Omega}} \|^2 \right) \nonumber \\
& \leq -R_{p-1}\left[(\gamma_p-\| \hat{A}_{p-1,p} \| )\sum_{i=1}^{r_p} \fourIdx{p}{}{2}{i}{\Omega}  \right.\nonumber \\
& \left. + (\gamma_{p-1}-\| \hat{A}_{p-1,p} \| ) \sum_{i=1}^{r_{p-1}}\fourIdx{p-1}{}{2}{i}{\Omega} \right]-\frac{R_{p-1,max}}{R_{p,min}}\gamma_p R_p\sum_{i=1}^{r_p} \fourIdx{p}{}{2}{i}{\Omega} \nonumber \\
& \leq -R_{p-1}\left[(\gamma_p-\| \hat{A}_{p-1,p} \| )\sum_{i=1}^{r_p} \fourIdx{p}{}{2}{i}{\Omega} \right. \nonumber \\
& \left. + \frac{(\gamma_{p-1}-\| \hat{A}_{p-1,p} \| )}{2} \sum_{i=1}^{r_{p-1}}\fourIdx{p-1}{}{2}{i}{\Omega} \right] \nonumber \\
& -\min\left\{ \frac{\gamma_p}{2},\frac{(\gamma_{p-1}-\| \hat{A}_{p-1,p} \| )}{2}\right\}V_{[p,p-1]}
\end{align}
where the second inequality has been arrived at by using available bounds on $R_p,R_{p-1}$ and applying the Cauchy-Schwarz inequality on the mixed term.

The following energy-like function allows amalgamation of the $p$ to $p-2$ blocks,
\begin{equation} \label{V13}
	V_{[p,p-2]} \triangleq V_{p-2}^o + \frac{R_{p-2,max}}{R_{p-1,min}} V_{[p,p-1]}
\end{equation} 
The pattern followed in prescribing the amalgamated energy function	is evident from the above equation. Proceeding along identical steps as Eq.~(\ref{V12d}), it can be shown that the directional derivative of $V_{[p,p-2]}$ turns out to be,
\begin{align} \label{V13d}
&	\dot{V}_{[p,p-2]} \leq -\frac{R_{p-2,max}}{R_{p-1,min}}\min\left\{ \frac{\gamma_p}{2},\frac{(\gamma_{p-1}-\| \hat{A}_{p-1,p} \| )}{2}\right\}V_{[p,p-1]}  \nonumber \\ 
& -R_{p-2}\Big[(\gamma_p-\| \hat{A}_{p-1,p} \|- \|\hat{A}_{p-2,p}\|)\|\fourIdx{p}{}{}{}{\bm{\Omega}} \|^2 \nonumber \\
& + \left(\frac{\gamma_{p-1}-\|\hat{A}_{p-1,p}\|}{2} - \|\hat{A}_{p-2,p-1}\|\right)\|\fourIdx{p-1}{}{}{}{\bm{\Omega}} \|^2 \nonumber \\
& + (\gamma_{p-2}-\|\hat{A}_{p-2,p-1}\|-\|\hat{A}_{p-2,p}\|)\|\fourIdx{p-2}{}{}{}{\bm{\Omega}} \|^2\Big] \nonumber \\
&\leq -\min\left\{ \frac{\gamma_p}{2},\frac{(\gamma_{p-2}-\|\hat{A}_{p-2,p-1}\|-\|\hat{A}_{p-2,p}\|)}{2},\right.  \nonumber \\ 
& \left. \frac{(\gamma_{p-1}-\| \hat{A}_{p-1,p} \| )}{2}\right\}V_{[p,p-2]} -R_{p-2}\Big[(\gamma_p-\| \hat{A}_{p-1,p} \| \nonumber \\
& - \|\hat{A}_{p-2,p}\|)\|\fourIdx{p}{}{}{}{\bm{\Omega}} \|^2 \nonumber \\
& + \bigg(\frac{\gamma_{p-1}-\|\hat{A}_{p-1,p}\|}{2} \|- \|\hat{A}_{p-2,p-1}\bigg)\|\fourIdx{p-1}{}{}{}{\bm{\Omega}} \|^2 \nonumber \\
&+ \frac{(\gamma_{p-2}-\|\hat{A}_{p-2,p-1}\|-\|\hat{A}_{p-2,p}\|)}{2}\|\fourIdx{p-2}{}{}{}{\bm{\Omega}} \|^2\Big]	
\end{align}										

Continuing in this prescribed manner the rest of the amalgamated Lyapunov candidate functions are defined by the following recursive formula,
\begin{equation} \label{V1j}
	V_{[p,p-j]} \triangleq V_{p-j}^o + \frac{R_{p-j,max}}{R_{p-j+1,min}} V_{[p,p-j+1]}
\end{equation} 
for $j=2,\ldots,p-1$. Diligently carrying out the derivatives of each $V_{[p,p-j]}$ and proceeding as before to compute $\dot{V}_{[p,p-j-1]}$, the following final candidate Lyapunov function can be arrived at,
\begin{equation} \label{miV}
	V \triangleq V_{[p,1]}
\end{equation} 
for which the directional derivative along dynamics (\ref{mizdyn}) can be compactly written as,
\begin{align} \label{miVd}
&	\dot{V} \leq  -\sigma V -R_1\left(\gamma_p - \sum_{l=1}^{p-1} \|\hat{A}_{p-l,p}\|\right)\|\fourIdx{p}{}{}{}{\bm{\Omega}}\|^2\nonumber \\
&				 -R_1\sum_{j=2}^{p-1}\left(\frac{\gamma_j}{2} - \frac{1}{2}\sum_{l=1}^{p-j} \|\hat{A}_{j,j+l}\|- \sum_{l=1}^{j-1} \|\hat{A}_{j-l,j}\|\right)\|\fourIdx{j}{}{}{}{\bm{\Omega}}\|^2 \nonumber \\
& 
\end{align}	
where $\sigma$ is defined in Eq.~(\ref{sig}). 
Assuming now that all inequalities~(\ref{gamjineq}) are satisfied, integrating both sides of Eq.~(\ref{miVd}) yields that $V(t) \leq V(0)\exp(-\sigma t)$. Further from the positive definiteness of each component function $V_j^o,\,j=1,2,\ldots,p$, it is possible to proceed backwards progressively starting at Eq.~(\ref{V1j}) to recover exponential convergence of each term $V_j^o$ at arbitrary rate $\sigma$. For example, 
\begin{equation} \label{conv1}
	V = V_{[p,1]} =  V_1^o + \big( R_{1,max}/R_{2,min} \big) V_{[p,2]}
\end{equation} 
which implies exponential convergence of $V_1^o$ and $V_{[p,2]}$ at rate $\sigma$. Then from the definition of $V_{[p,2]}$, i.e.,
\begin{equation} \label{conv2}
	 V_{[p,2]} =  V_2^o + \big(R_{2,max}/R_{3,min}\big) V_{[p,3]}
\end{equation} 
exponential convergence of $V_2^o$ and $V_{[p,3]}$ at an identical rate can be concluded. Similarly, subsequent steps will prove exponential convergence with rate $\sigma$ for all $V_j^o$. This along with the fact that there exists an $R_{j,min} > 0$ corresponding to each $R_j(t)$ implies exponential convergence of states $\fourIdx{j}{}{}{}{\bm{\Omega}} = (\sqrt{V_j^o}/R_{j,min})$ at a rate $\sigma/2$. The invertibility of the augmented states definitions (\ref{Om1}), (\ref{Om2}) and (\ref{Omj}) to recover $\fourIdx{j}{}{}{}{\bm{z}}$ further proves exponential convergence of $\fourIdx{j}{}{}{}{\bm{z}}$ and in turn that of $\xv$ to zero at the same rate.
\end{proof}
%

\subsection{Adaptive Control Modification}

It is possible to introduce adaptation to the control law described in Theorem~\ref{miPEthm1} for situations wherein the matrices describing the system, i.e. $A$ and $B$ are unknown.  A controllability assumption is made on the matrix pair $(A,B)$ as in Theorem~\ref{miPEthm1}. Further, these matrices are assumed to be unknown but in the canonical form described by Eq.~(\ref{mizdyn}). The following Corollary is therefore stated in terms of matrices $\hat{A}$ and $\hat{B}$ as system matrices. An important facet of this assumption is that the minimum number of controls required to ascertain system controllability is known (this is the same as the number of blocks in $\hat{A}$ denoted by $p$).

To begin with, a new persistence filter state $R(t) \in \mathbb{R}$ is proposed with bounded, time-varying gains and lower boundedness of $R(t)$ as in Lemma~4 of \citep{srikant09} for the constant $\lambda$ case shown.
\begin{lem} \label{rmint}
Consider the persistence filter defined by,
\begin{equation} \label{miRt}
    \dot{R} = -\hat{\lambda}(t) R + g^k(t), 
\end{equation}
	Let $g(\cdot):\mathbb{R} \rightarrow \mathbb{R}$ be $\mathcal{C}^{n-1}$, bounded with bounded derivatives up to order $(n-1)$, $n$ being the order of the dynamics and PE, $k \triangleq \textrm{max}\{2,2^{\lceil \log_2 n\rceil}\}$, $0 < \hat{\lambda}_{min} \leq \hat{\lambda}(t) \leq \hat{\lambda}_{max} $. Then the solution $R(\cdot)$ of (\ref{miRt}) with initial condition $R(0) \, > \, 0$ satisfies,
	\begin{equation} \label{rmt}
		\exists R_{min},R_{max}>0 \,\,s.t. \,\, \forall t \geq 0 : \,\,0<R_{min} \leq R(t) \leq R_{max}
	\end{equation}
\end{lem}

The following corollary to Theorem~\ref{miPEthm1} stated without proof formalizes the adaptive control extension.
\begin{cor} \label{miPEthm2}
Consider the linear, multi-input dynamics in Eq.~(\ref{mizdyn}) with unknown system matrices $\hat{A}$ and $\hat{B}$ assumed to form a controllable pair and the component functions of $G(t)$ are $\mathcal{C}^{n-1}$, bounded with bounded derivatives up to order $(n-1)$ and satisfy the PE condition. Further, let the augmented states be defined as in Eqs.~(\ref{Om1}) and (\ref{Omj}). Then the following control law,
\begin{equation} \label{mipcont2}
        u_j = -\frac{1}{g_j}\sum_{i=1}^{r_j}\hat{\alpha}_{j,r_j-i+1}(\fourIdx{j}{}{}{i}{\Omega}-\fourIdx{j}{}{}{i}{z})- \frac{1}{g_j}\frac{\ud}{\ud t} (\fourIdx{j}{}{}{r_j}{\Omega}-\fourIdx{j}{}{}{r_j}{z}) - \frac{g^{k-1}_j}{2R_j}\fourIdx{j}{}{}{r_j}{\Omega}	
\end{equation}
for $j=1,2,\ldots,p$ and $u_j = 0$ for $j>p$ with the persistence filters defined by,
\begin{equation} \label{miRtj}
    \dot{R}_j = -\hat{\lambda}_j(t) R_j + g_j^k(t), 
\end{equation}
with $R_j(0) > 0$, $k = \textrm{max}\{2,2^{\lceil \log_2 r \rceil}\}$, $r = \text{max} \{r_1,r_2,\ldots,r_p\}$, the state $\hat{\lambda}_j(t)$ has dynamics,
\begin{equation} \label{lamhatmi}
	\dot{\hat{\lambda}}_j = \nu_j R_j \| \fourIdx{j}{}{}{}{\bm{\Omega}} \|^2, \quad \text{ any } \nu_j,\hat{\lambda}_j(0) > 0
\end{equation}	
and parameter estimates evolving as,
\begin{align} \label{ajhat}
&	\dot{\hat{\alpha}}_{j,r_j-i+1} = - 2\eta_{j,i}R_j\fourIdx{j}{}{}{r_j}\Omega\sum_{i=1}^{r_j} \tilde{\alpha}_{j,r_j-i+1}(\fourIdx{j}{}{}{i}{\Omega}-\fourIdx{j}{}{}{i}{z}), \nonumber \\
& \text{ any } \hat{\alpha}_{j,r_j-i+1}(0),,\eta_{j,i}>0 \quad j = 1,2,\ldots,p
\end{align}
guarantees asymptotic convergence of $\zv$ to the origin.
\end{cor}

\section{Some Remarks on Observer Design for Multi-Output Systems} \label{mobs}


The successful resolution of the multi-input stabilization problem with a diagonal singular gain matrix by transformation to a corresponding block-triangular canonical form allows for a similar observer counterpart for multi-output systems. To this end the following input-output dynamics is considered,
\begin{eqnarray} \label{miodyn}
	\dot{\xv} &=& A\xv + \bm{u} \nonumber \\
	y &=& G(t)C\xv
\end{eqnarray}
with $\xv \, \in \mathbb{R}^n$, $\bm{u} \in \mathbb{R}^n$, $A \in \mathbb{R}^{n\times n}$, $C \in \mathbb{R}^{s \times n}$ and $G(\cdot):\mathbb{R}^{\geq 0} \rightarrow \mathbb{R}^{s \times s}$. It is assumed that the $(A,C)$ pair is observable and $G(t)$ is sufficiently smooth, satisfies the PE condition~\citep[p. 24-25]{sastry} and has the form,
\begin{equation} \label{Go}
	G(t) = \left[\begin{array}{cccc}
								g_1(t) & 0 & \ldots & 0 \\
								0 & g_2(t) & \ldots & 0 \\
								\vdots & \vdots & \ddots & \vdots \\
								0 & 0 & \ldots & g_s(t)
								\end{array} \right]		
\end{equation}	
The PE claim on $G(t)$ implies as before that each of the component functions $g_i(t),\,i=1,2,\ldots,s$ satisfy the PE condition. The observer design for single-output dynamics with singular measurement gains was developed in \citet{srikant12j}. Here we extend it to the multi-output case.

Observability of the $(A,C)$ pair in Eq.~(\ref{miodyn}) implies controllability of $(A^T,C^T)$. The canonical transformation described in Section~\ref{PEstab} can therefore be constructed to obtain a non-singular matrix $T$ such that, $\hat{A} = TA^TT^{-1}$ and $\hat{B} = TC^T$ where $\hat{A}$ and $\hat{B}$ are defined in Eqs.~(\ref{Ahat})-(\ref{Bhat}). The similarity transformation, $\bm{z} = (T^{-1})^T \xv$ is now considered. This results in the following transformed dynamics,
\begin{eqnarray} \label{mizodyn}
		\dot{\zv} &=& A_o\zv + (T^{-1})^T \bm{u} \nonumber \\
		y &=& G(t)C_o\zv
\end{eqnarray}
with $A_o = \hat{A}^T$ and $C_o = \hat{B}^T$. It has therefore been possible to obtain canonical dynamics which are the transpose of the control problem case. This is now in suitable form to carry out observer design. Following are some of the features of the transformed system,
\begin{itemize}

\item The matrix $A_o$ is now upper triangular so that the $1$-block is decoupled as opposed to the $p$-block in the control case. 

\item Each diagonal block has the same structure as the observer canonical form for a single output system with the last state in each corresponding block being the measured output. More specifically the first $p$ outputs are, $y(1:p) = [g_1\fourIdx{1}{}{}{r_1}{z},\: g_2\fourIdx{2}{}{}{r_2}{z},\: \ldots,\: g_p\fourIdx{p}{}{}{r_p}{z}]^T$ assuming that $\zv$ is still partitioned according to Eq.~(\ref{zpart}).

\item The last $s-p$ rows in the $C_o$ are redundant outputs and need not be used in the observer design process since the system is observable with only the first $p$ outputs.	

\item The upper triangular structure of $A_o$ implies that each block $(1+j)$-block has coupling with the $j$-block for $j=1,2,\ldots,p-1$. However based on the structure of $\hat{A}_{k,j}$ in Eq.~(\ref{Akj}) it can be inferred that the coupling terms appear only in the first state corresponding to each block.

\item The evident similarities with the dual control design problem leads to the belief that an exponential observer design for multi-input systems with diagonal, singular gains on the measurements can be accomplished by extending Theorem~1 of \citep{srikant12j} based on the above canonical form.

\item The observer design would begin at the uncoupled $1$-block estimation error dynamics and then proceeding in steps to higher numbered blocks which is exactly the reverse of the order followed in the control design. The estimation error dynamics corresponding to the $2$-block has coupling with the error states in the $1$-block which have already been rendered exponentially stable. This implies that the coupling terms wash out exponentially fast resulting in exponential convergence of the $2$-block error states. These arguments can potentially be continued to establish exponential reconstruction of all states.

\end{itemize}

\subsection{Observer Design}

Let the observer have the following dynamics,
\begin{eqnarray} \label{mizobs}
	\dot{\hat{\zv}} &=& A_o\hat{\zv} + (T^{-1})^T \bm{u} + \bm{\eta} \nonumber \\
	\hat{y} &=& G(t)C_o\hat{\zv}
\end{eqnarray}
Then the error dynamics for $\tilde{\zv} \triangleq \zv - \hat{\zv}$ is,
\begin{eqnarray} \label{mizerr}
	\dot{\tilde{\zv}} &=& A_o\tilde{\zv} - \bm{\eta} \nonumber \\
	\tilde{y} &=& G(t)C_o\tilde{\zv}
\end{eqnarray}	
with $\bm{\eta} = [\fourIdx{p}{}{}{}{\bm{\eta}},\: \fourIdx{p-1}{}{}{}{\bm{\eta}},\cdots,\:\fourIdx{1}{}{}{}{\bm{\eta}}]^T$.

The analysis begins by looking at the error dynamics of the $1$-block which as stated earlier is an uncoupled single-output system and can be verified to have the following form,
\begin{eqnarray} \label{mizerr1}
	\frac{\ud}{\ud t}\fourIdx{1}{}{}{}{\tilde{\zv}} &=& \hat{A}^T_{1,1}\,\fourIdx{1}{}{}{}{\tilde{\zv}} - \fourIdx{1}{}{}{}{\bm{\eta}} \nonumber \\
	\tilde{y}_1 &=& g_1 \,\fourIdx{1}{}{}{r_1}{z}
\end{eqnarray}
Therefore applying Theorem~1 of \citep{srikant12j} to design $\fourIdx{1}{}{}{}{\bm{\eta}}$ guarantees that $\fourIdx{1}{}{}{}{\zv} \rightarrow \bm{0}$ exponentially.

Going forward, the error dynamics of the $2$-block has a coupling with $\fourIdx{1}{}{}{}{\tilde{\zv}}$ as follows,
\begin{eqnarray} \label{mizerr2}
	\frac{\ud}{\ud t}\fourIdx{2}{}{}{}{\tilde{\zv}} &=& \hat{A}^T_{2,2}\,\fourIdx{2}{}{}{}{\tilde{\zv}} + \hat{A}^T_{1,2}\,\fourIdx{1}{}{}{}{\tilde{\zv}} - \fourIdx{2}{}{}{}{\bm{\eta}} \nonumber \\
	\tilde{y}_2 &=& g_2 \,\fourIdx{2}{}{}{r_2}{z}
\end{eqnarray}	
Suppose $\fourIdx{2}{}{}{}{\bm{\eta}}$ is chosen based on Theorem~1 of \citep{srikant12j} for the following nominal dynamics (i.e., if the coupling term did not exist),
\begin{eqnarray} \label{mizerr2nom}
	\frac{\ud}{\ud t}\fourIdx{2}{}{}{}{\tilde{\zv}} &=& \hat{A}^T_{2,2}\,\fourIdx{2}{}{}{}{\tilde{\zv}} - \fourIdx{2}{}{}{}{\bm{\eta}} \nonumber \\
	\tilde{y}_2 &=& g_2 \,\fourIdx{2}{}{}{r_2}{z}
\end{eqnarray}	
Substituting, the innovation term $\fourIdx{2}{}{}{}{\bm{\eta}}$ thus obtained back in Eq.~(\ref{mizerr2}) it can be concluded that in the absence of the coupling term the linear, time-varying system is exponentially stable. Further, it is known from the analysis corresponding to the $1$-block that $\fourIdx{1}{}{}{}{\tilde{\zv}}$ and hence the coupling term is exponentially decaying. Therefore, the observer dynamics is an exponentially stable linear system forced by an exponentially decaying function, thus allowing the conclusion that $\fourIdx{2}{}{}{}{\tilde{\zv}} \rightarrow 0$ exponentially.

The analysis can be continued along identical lines to conclude exponential convergence of all $\fourIdx{i}{}{}{}{\tilde{\zv}}$ for $i=1,2,\ldots,p$ to the origin. This leads to the following Lemma on observer design for multi-output dynamics with matrix, singular, time-varying gains.

\begin{lem} \label{miobs}
	Consider the linear, multi-output dynamics in Eq.~(\ref{miodyn})-(\ref{Go}) under the assumptions that $(A,C)$ is an observable pair and the component functions of $G(t)$ are $\mathcal{C}^{n-1}$, bounded with bounded derivatives up to order $(n-1)$ and satisfy the PE condition. Further, let $\bm{z} = (T^{-1})^T\xv$ transform the dynamics~(\ref{miodyn}) to Eq.~(\ref{mizodyn}). Then the following observer,
	\begin{equation} \label{mixobs}
		\dot{\hat{\xv}} = A\hat{\xv} + \bm{u} + T^T \bm{\eta}
	\end{equation}
with $\bm{\eta} = [\fourIdx{p}{}{}{}{\bm{\eta}},\: \fourIdx{p-1}{}{}{}{\bm{\eta}},\cdots,\:\fourIdx{1}{}{}{}{\bm{\eta}}]^T$ where each $\fourIdx{i}{}{}{}{\bm{\eta}}$ ($i=1,2,\ldots,p$) is designed using Theorem~1 of \citep{srikant12j} on the following nominal dynamics,
	\begin{eqnarray} \label{mizerrinom}
	\frac{\ud}{\ud t}\fourIdx{i}{}{}{}{\tilde{\zv}} &=& \hat{A}^T_{i,i}\,\fourIdx{i}{}{}{}{\tilde{\zv}} - \fourIdx{i}{}{}{}{\bm{\eta}} \nonumber \\
	\tilde{y}_i &=& g_i \,\fourIdx{i}{}{}{r_i}{z}
\end{eqnarray}
guarantees $\tilde{\xv} \triangleq (\xv - \hat{\xv}) \rightarrow \bm{0}$ exponentially under the assumptions made in the aforementioned Theorem.
\end{lem}
	
The proposed observer is of lower order than the corresponding Kalman Filter based estimator. However, it does not share similar optimality guarantees in presence of noise. A notable advantage of the proposed observer is that there is much greater control over the exponential convergence rate (dictated by the persistence filter bandwidth) as compared to Kalman Filtering based techniques. The primary target of this observer construction is to explore nonlinear extensions of the persistence filter based observer design. 

\section{Axi-Symmetric Underactuated Spacecraft Stabilization} \label{sim}

An example of application of Theorems~\ref{miPEthm1} and Corollary~\ref{miPEthm2} to stabilization of a nonlinear dynamical system is considered in this section. This illustrates the possibility of extension of the results stated in previous sections to special nonlinear dynamical systems. Here we consider the attitude stabilization of a spacecraft with a single axis of symmetry and only two independent actuators. 

It is assumed without loss of generality that the body axis of the spacecraft is aligned with the principal axis and further that there is an axis of symmetry, say $J_2 = J_3$ where $J_1, J_2, J_3$ are the principal moments of inertia. It is assumed that there are only two physical actuation mechanisms one of which can however reorient to alternately provide torque on two axes (ref. `Thruster Gimballing' in \citet{stanton09}). The linearized attitude kinematics and angular velocity dynamics for this setup written in the body frame of reference are,
\begin{eqnarray} 
	\dot{q}_0 &=& 0 \label{att0}\\
	\dot{\bm{q}}_v &=& \frac{1}{2}\bm{\omega} \label{attv} \\
	\dot{\omega}_1 &=& \frac{g_1(t)}{J_1} u_1 \label{angv1} \\
	\dot{\omega_2} &=& \frac{(J_3-J_1)}{J_2} \omega_1\omega_3 + \frac{g_2(t)}{J_2} u_2 \label{angv2} \\
	\dot{\omega_3} &=& \frac{(J_1-J_2)}{J_3} \omega_1\omega_2 + \frac{g_2(t)}{J_3} u_3 \label{angv3}
\end{eqnarray}
where the attitude is represented using Euler parameters $[q_0,\,\bm{q}_v^T]^T$ ~\citep{schaub03}. The kinematics equations (\ref{att0})-(\ref{attv}) have been linearized while the angular velocity dynamics (\ref{angv1})-(\ref{angv3}) represent the complete nonlinear dynamics. The dynamics corresponding to $\omega_2$ and $\omega_3$ have the same gain on the control representing identical actuator schedules. This assumption can be relaxed and the gain on $\omega_3$ dynamics can be 	set to unity representing a dedicated actuator in the subsequent analysis without effecting the results. Further, the gains $g_1(t)$ and $g_2(t)$ are designed to be time-wise orthogonal functions to signify actuator reorientation. It can be assumed without loss of generality that one of the actuators can reorient to provide torques about the first and second principal axes. The other actuator is assumed to provide torque about the third principal axis but has the same firing schedule as in the second axis ($g_2(t)$). An example could be a cylindrical satellite with a reorienting actuator (say a gimballing thruster) between the roll and yaw directions while having a fixed actuator in the pitch direction. Reducing the number of actuators required for angular velocity stabilization has obvious advantages especially with reference to micro-satellites with considerable weight and power constraints.

The stabilization objective can be formalized in terms of state-variables in Eq.~(\ref{attv})-(\ref{angv3}) as $\bm{q}_v,\,\bm{\omega} \rightarrow 0$ as $t\rightarrow \infty$ where $\bm{\omega} = [\omega_1,\, \omega_2,\,\omega_3]^T$. The sub-system
\begin{displaymath}
	\dot{q}_{v1} = 0.5\, \omega_1; \; \dot{\omega}_1 = g_1(t) u_1 /J_1
\end{displaymath}
of the dynamics Eq.~(\ref{attv})-(\ref{angv1}) is an uncoupled single-input system and so direct application of Theorem~8 of \citep{srikant09} allows design of control law $v_1 \triangleq u_1 /J_1$ and persistence filter state $R_1$ that guarantees that $|\omega_1(t)| \leq \alpha \exp{-(\beta t)}$. In this convergence estimate although the rate $\beta$ is known and can be chosen arbitrarily large, the knowledge of $R_{1,min}$ and hence $\alpha$ (which depends on $R_{1,min}$) is not known~\citep{srikant09}. Defining placeholders, $k_2 = (J_3-J_1)/J_2$, $k_3 = (J_1-J_2)/J_3$ and modified controls, $v_2 = u_2/J_2$ and $v_3 = u_3/J_3$, the remaining dynamics can be re-written using Eq.~(\ref{attv})-(\ref{angv3}) as,
\begin{eqnarray} 
	\dot{q}_{v2} &=& \frac{1}{2}\omega_2 \label{attv2} \\
	\dot{\omega_2} &=& k_2 \omega_1\omega_3 + g_2(t)v_2 \label{angv2n} \\
	\dot{q}_{v3} &=& \frac{1}{2}\omega_3 \label{attv3} \\
	\dot{\omega_3} &=& k_3 \omega_1\omega_2 + g_2(t)v_3 \label{angv3n}
\end{eqnarray}

The following nonlinear persistence filter is now defined along the lines of Eq.~(\ref{miRtj}),
\begin{equation} \label{angR2}
	\dot{R}_2 = -\hat{\lambda}_2 R_2 + g_2^2(t)
\end{equation}
with $R_2(0),\,\hat{\lambda}_2(0) > 0$. We now define $V_2 \triangleq R_2 [q_{v2}^2 + \Omega_2^2]$ with augmented state $\Omega_2$ defined identical to the single-input case as,
\begin{displaymath}
	\Omega_2 = \omega_2 + \frac{g_2^2}{R_2} q_{v2}
\end{displaymath}
	
Then with the following control law defined along the same lines as the single-input case,
\begin{equation} \label{v2}
	v_2 = -\frac{g_2}{2R_2} \Omega_2 + k_2\frac{g_2}{R_2}\omega_1q_{v3} - \frac{1}{g_2} \frac{\ud}{\ud t} (\Omega_2 - \omega_2)
\end{equation}
yields,
\begin{eqnarray}
	\dot{V}_2 &=& -\hat{\lambda}_2 R_2(\Omega_2^2 + q_{v2}^2) + R_2q_{v2}\Omega_2 + 2k_2R_2 \omega_1\Omega_2\Omega_3 \\
						&\leq & -(\hat{\lambda}_2 - 0.5)R_2q_{v2}^2 - \left(\hat{\lambda}_2 - k_2\alpha - 0.5\right)R_2 \Omega_2^2 + k_2\alpha\Omega_3^2 \nonumber \\
						& &
\end{eqnarray}
where the bound on $|\omega_1|$ and the inequalities (i) $\alpha \exp(-\beta t) \leq \alpha$, (ii) $2\|a\| \|b\| \leq \|a\|^2 + \|b\|^2$ have been used to arrive at the second inequality.

Similarly, with a choice of $V_3 \triangleq R_2 [q_{v3}^2 + \Omega_3^2]$ and augmented state $\Omega_3$ defined as,
\begin{displaymath}
	\Omega_3 = \omega_3 + \frac{g_2^2}{R_2} q_{v3}
\end{displaymath}
and control,
\begin{equation} \label{v3}
	v_3 = -\frac{g_2}{2R_2} \Omega_3 + k_3\frac{g_2}{R_2}\omega_1q_{v2} - \frac{1}{g_2} \frac{\ud}{\ud t} (\Omega_3 - \omega_3)
\end{equation}
results in,
\begin{equation}
	\dot{V}_3 \leq  -(\hat{\lambda}_2 - 0.5)R_2q_{v3}^2 - \left(\hat{\lambda}_2 - k_3\alpha - 0.5\right)R_2 \Omega_3^2 + k_3\alpha\Omega_2^2
\end{equation}

Then by choosing the consolidated Lyapunov candidate as, $V \triangleq V_1 + V_2 + \tilde{\lambda}_2^2/2\gamma$ with $\tilde{\lambda}_2 = \hat{\lambda}_2 - k_2\alpha - k_3 \alpha - 0.5 - \epsilon$ for some $\epsilon > 0$, the derivative along the closed-loop trajectories of Eq.~(\ref{attv2})-(\ref{v3}) satisfies,
\begin{align}
&	\dot{V} \leq -(\hat{\lambda}_2 - k_2\alpha - k_3 \alpha - 0.5)R_2 (\Omega_2^2 + \Omega_3^2)\nonumber \\
&  - (\hat{\lambda}_2 - 0.5)R_2(q_{v2}^2 + q_{v3}^2) + \frac{1}{\gamma}\tilde{\lambda}_2 \dot{\hat{\lambda}}_2 \\
& \leq -(\hat{\lambda}_2 - k_2\alpha - k_3 \alpha - 0.5) R_2 (q_{v2}^2 + q_{v3}^2 + \Omega_2^2 + \Omega_3^2)  + \frac{1}{\gamma}\tilde{\lambda}_2 \dot{\hat{\lambda}}_2
\end{align}					

If we choose,
\begin{equation}
	\dot{\hat{\lambda}}_2 = \gamma R_2 (q_{v2}^2 + q_{v3}^2 + \Omega_2^2 + \Omega_3^2)
\end{equation} 
similar to Eq.~(\ref{lamhatmi}) we get, $	\dot{V}  \leq -\epsilon R_2(q_{v2}^2 + q_{v3}^2)$. Beyond this, application of standard signal-chasing and Barbalat's Lemma arguments along with the fact that $R_2$ is lower bounded above zero implies asymptotic convergence of $q_{v2},\,q_{v3},\,\omega_2,\,\omega_3$ to zero. This along with exponential convergence of $q_{v1},\,\omega_1$ to zero from before completes the proof.		

\begin{figure}
		\centering
		\includegraphics[width = 0.7\textwidth, keepaspectratio]{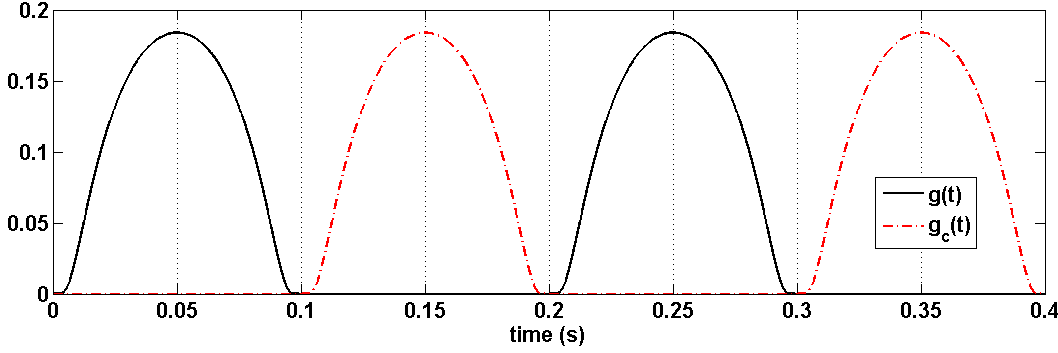}
	\caption{Conjugate smooth functions $g_1(t)$ and $g_2(t)$}\label{ggc}
\end{figure}

A sample simulation on a representative micro-satellite was carried out to test the efficacy of the proposed control. For the purposes of the simulation the inertias were chosen as $J_1 = 3 \:kgm^2$, $J_2=J_3 = 2 \:kgm^2$. The initial attitude of the satellite is a rotation of 18$^\circ$ about the axis $[1/\sqrt{3},\: 1/\sqrt{3},\: 1/\sqrt{3}]^T$, the initial angular velocities are $\bm{\omega}(0) = 0.1[\pi/12,\: -\pi/6, \: \pi/8]^T$ and other states initialized at $R_1(0) = R_2(0) = 1$, $\hat{\lambda}_2(0) = 2$. The various simulation parameters were chosen as $\lambda_1 = 2$ and $\nu = 0.01$. Further, the scheduling functions $g_1(t)$ and $g_2(t)$ were chosen orthogonal to each other with $g_1(t)$ active over 1.8 s followed by a 0.4 s gap where none of the actuators are operating followed by a 1.8 s $g_2(t)$ `on' window over a cycle of 4 s. This accounts for actuator reorientation as required. The design of $g_1(t)$ and $g_2(t)$ is based on smooth, compactly supported functions described in \citet{jamshidi06} to define the conjugate functions $g_1(t)$ and $g_2(t)$,
\begin{equation} \label{gi}
	\phi(t) = (1+ \cos(t\pi))H(1-t^2)
\end{equation}	
and look similar in nature to Fig.~\ref{ggc}.

\begin{figure}
		\centering
		\includegraphics[width = 0.8\textwidth, keepaspectratio]{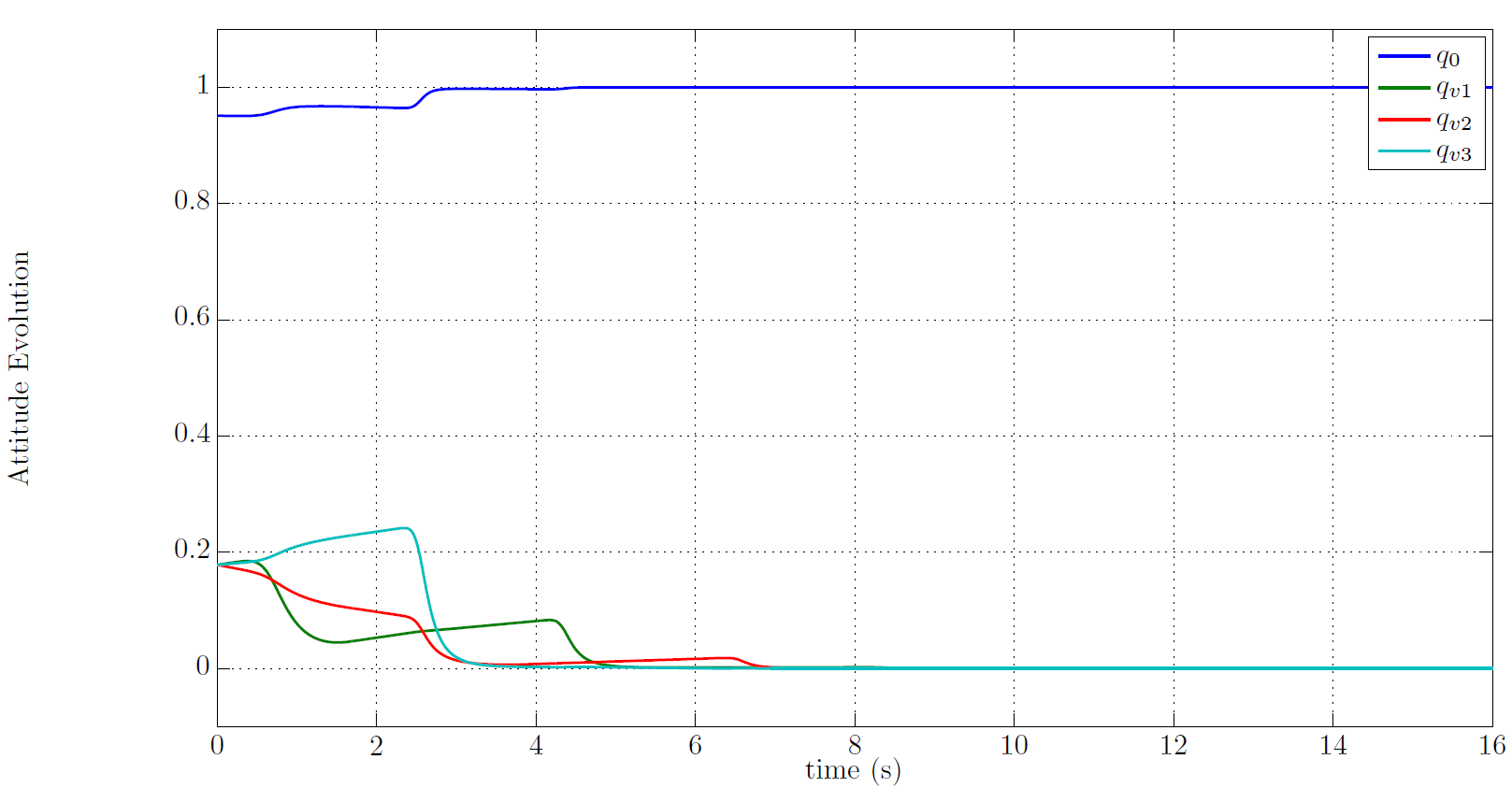}
	\caption{Quaternions - Attitude Stabilization of an Underactuated Spacecraft.}\label{qnion}
\end{figure}

\begin{figure}
		\centering
		\includegraphics[width = 0.8\textwidth, keepaspectratio]{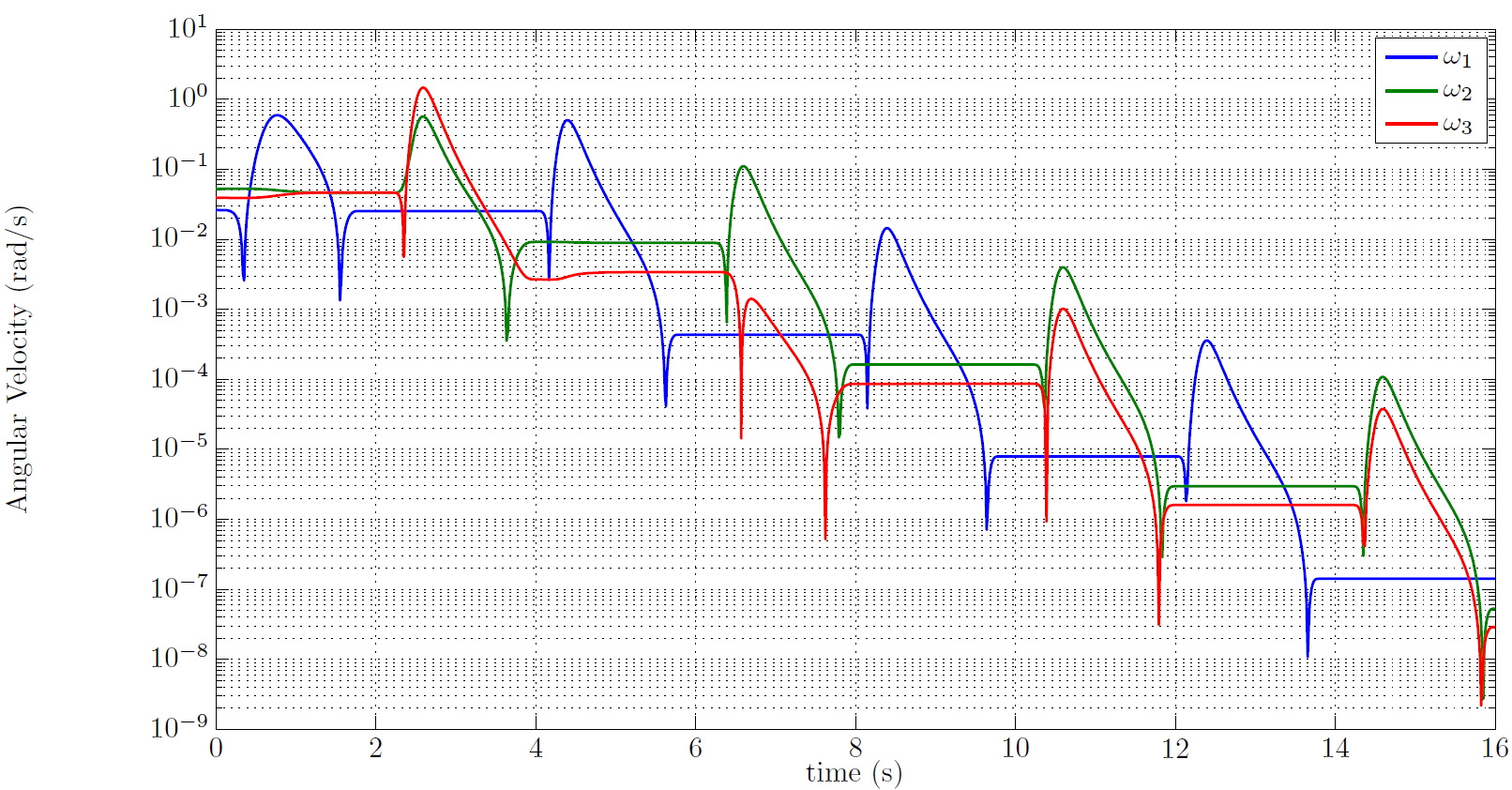}
	\caption{Angular Velocity - Attitude Stabilization of an Underactuated Spacecraft.}\label{omega}
\end{figure}

\begin{figure}
		\centering
		\includegraphics[width = 0.8\textwidth, keepaspectratio]{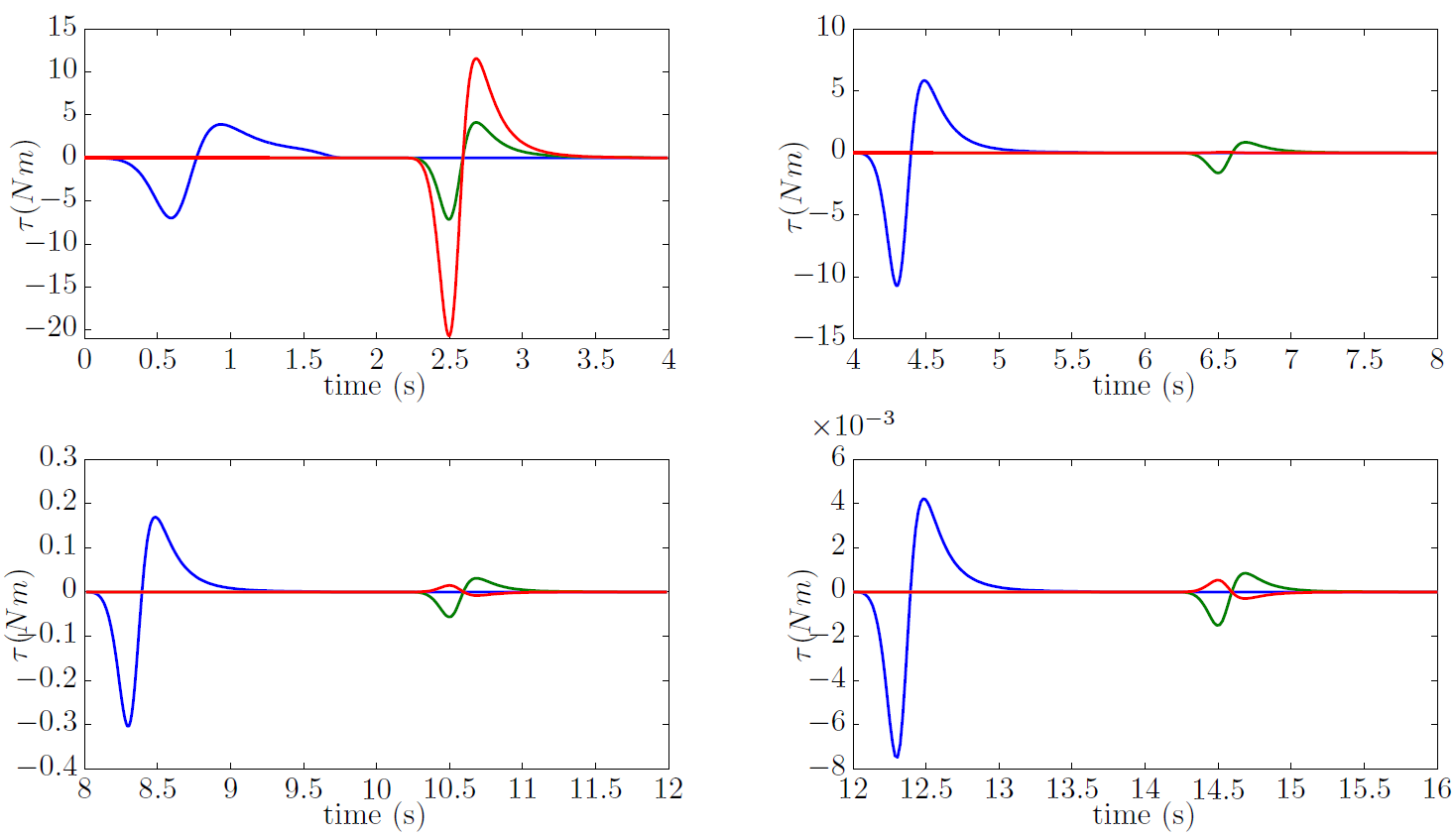}
	\caption{Control Torques - Angular Velocity Stabilization for Underactuated Spacecraft. Blue, Green and Red curves represent control commands on first, second and third principal axes respectively.}\label{ontrols}
\end{figure}

Although the control design is based on linearized kinematics~(Eq.~\ref{attv}), the simulations apply the control law to the true nonlinear attitude kinematics $(\dot{q}_0 = -0.5*q_v^T\omega;\,\dot{q}_v = 0.5*q_0\omega + 0.5*q_v\times \omega)$ and dynamics equations. Figure~\ref{qnion} shows the evolution of the attitude quaternions, Figure~\ref{omega} shows the angular velocity evolution and Figure~\ref{ontrols} the corresponding control torque in each actuation cycle. It is evident from these plots that the control law produces the desired convergence of the quaternions and angular velocities with reasonable control torques and using only two reorientable actuators. Further, it is evident from Figure~\ref{ontrols} that one of the actuators alternatively actuates both the first and the second body axis of the satellite. 

\section{Concluding Remarks} \label{conc}
This article shows extension of the persistence filter based controller framework to stabilize multi-input, linear dynamics with a diagonal, time-varying control gain. The number of controls is allowed to be less than the number of states as long as the system satisfies the linear time invariant controllability condition in absence of the time-dependent control gains. The gains can potentially pass through singular phases representing gaps in actuation. However, these time-varying gains are assumed to satisfy the persistence of excitation condition to allow for the system to be controllable in presence of the gain matrix. The stabilizing controller is designed by transforming the original multi-input system to a canonical form consisting of a series of single-input dynamics with uni-directional coupling. Subsequently, a set of persistence filters are defined corresponding to each single-input system in the canonical form. The structure of the control law is similar to the single-input case. An adaptive control result for the case wherein the system matrices are unknown was also developed. It was shown that a modified nonlinear persistence filter allows construction of stabilizing multi-input control law even when the plant matrices are not precisely known. The modified persistence filter formulation was employed to design a feedback law for asymptotically stabilizing the attitude and angular velocity of a spacecraft with only two actuators, one of which can reorient to alternately provide torque in two directions. The application demonstrates possibility of applying the linear, multi-input persistence filter based control design to nonlinear dynamical systems. The single-output observer design with singular, time-varying measurement gains has also been extended to the multi-output, special case of diagonal gains in this paper. This development uses an identical canonical transformation as the dual control counterpart. In future, the authors will look at extending the observer design result to nonlinear dynamics.

\section*{Acknowledgment}

This research work was supported in part by National Aeronautics and Space Administration Grant NNX09AW25G with Dr.~Timothy Crain as Program Manager.



\bibliographystyle{plainnat}
\bibliography{diss}
%

\end{document}